\documentclass[11pt,reqno]{amsart}
\usepackage{geometry}                
\usepackage{amsmath}
\usepackage{amssymb}
\usepackage{amsthm}
\usepackage{amsxtra}
\usepackage{amscd, amsrefs}
\usepackage{mathtools}
\usepackage[all]{xy}
\usepackage{blindtext}
\usepackage[utf8]{inputenc}
\usepackage[english]{babel}

\usepackage{tikz}

\usepackage{microtype}
\usepackage{esint}

 \usepackage{hyperref}
\DeclareMathOperator{\supp}{supp}

\theoremstyle{plain}
\newtheorem{thm}{Theorem}[section]
\newtheorem{lem}[thm]{Lemma}
\newtheorem{prop}[thm]{Proposition}
\newtheorem{cor}[thm]{Corollary}

\theoremstyle{definition}
\newtheorem{defn}[thm]{Definition}

\newtheorem{condition}[thm]{Condition}
\theoremstyle{remark}

\numberwithin{equation}{section}

\title[]{A global space-time estimate for dispersive operators through its local estimate}

\author[C.Cho]{Chu-Hee Cho}
\address{Department of Mathematical Sciences and RIM, Seoul National University, Seoul 151--747, Republic of Korea}
\email{akilus@snu.ac.kr}

\author[Y. Koh]{Youngwoo Koh}
\address{Department of Mathematics Education, Kongju National University, Kongju 32588, Republic of Korea}
\email{ywkoh@kongju.ac.kr}

\author[J. Lee]{Jungjin Lee}
\address{Department of Mathematical Sciences, School of Natural Science, Ulsan National Institute of Science and Technology, UNIST-gil 50, Ulsan 44919, Republic of Korea}
\email{jungjinlee@unist.ac.kr}
\subjclass{35A24, 35L10, 35Q41, 35S35}
\keywords{Schr\"odinger equation, dispersive operator, global space-time estimate}

\thanks{C. Cho is supported by NRF grant no.2020R1I1A1A01072942 (Republic of Korea). Y. Koh is supported by the research grant of the Kongju National University in 2021. J. Lee is supported by NRF grant no.2020R1I1A1A01062209 (Republic of Korea).
}
\date{\today}
\begin{document}

\begin{abstract}
We will show that a local space-time estimate implies a global space-time estimate for dispersive operators. In order for this implication we consider a Littlewood-Paley type square function estimate for dispersive operators in a time variable and a generalization of Tao's epsilon removal lemma in mixed norms. By applying this implication to the fractional Schr\"odinger equation in $\mathbb R^{2+1}$ we obtain the sharp global space-time estimates with optimal regularity from the previous known local ones. 
\end{abstract}

\maketitle

\section{Intoduction}
Let us consider a Cauchy problem of a dispersive equation in $\mathbb R^{n+1}$
    \begin{equation} \label{gen_form}
    \left\{
    \begin{aligned}
    i \partial_t u +  \Phi(D)u &= 0,\\ 
    u(0) &= f, 
    \end{aligned}
    \right.
    \end{equation}
where $\Phi(D)$ is the corresponding Fourier multiplier to the function $\Phi$. We assume that $\Phi \in C^\infty(\mathbb R^n \setminus \{0\})$ is a real-valued  function satisfying the following conditions:
\begin{condition} \label{homog_cond} 
{$ $}
\begin{itemize}
\item
$|\nabla \Phi(\xi)| \neq 0$ for all $\xi \neq 0$.
\item
There is a constant $\mu \ge 1$ such that $\mu^{-1} \le |\Phi(\xi)| \le \mu$ for any $\xi$ with $|\xi|=1$.
\item
There is a constant $m \ge 1$ such that $\Phi(\lambda\xi)=\lambda^m\Phi(\xi)$ for all $\lambda>0$ and all $\xi \neq 0$.
\item
The Hessian $H_{\Phi}(\xi)$ of $\Phi$ has rank at least $1$ for all $\xi \neq 0$.
\end{itemize}
\end{condition}
The solution $u$ to \eqref{gen_form} becomes the Schr\"odinger operator $e^{-it\Delta}f$ if $\Phi(\xi)=|\xi|^2$ and the wave operator $e^{it\sqrt{-\Delta}}f$ if $\Phi(\xi)=|\xi|$. When $\Phi(\xi)=|\xi|^m$ for $m>1$, the solution is called the fractional Schr\"odinger operator $e^{it(\sqrt{-\Delta})^{m/2}}f$. 

Let $e^{it\Phi(D)}f$ denote the solution to \eqref{gen_form}.  
Our interest  is to find suitable pairs $(q,r)$ which satisfy the global space-time estimate
    \begin{equation}\label{main_goal}
    \|e^{it\Phi(D)} f\|_{L_x^q(\mathbb R^n; L_t^r(\mathbb R))}
    \le C \|f\|_{\dot H^s(\mathbb R^n)},
    \end{equation}
where $\dot H^s(\mathbb R^n)$ denotes the homogeneous $L^2$ Sobolev space of order $s$.
By scaling invariance the regularity $s=s(r,q)$ should be defined as
\begin{equation} \label{regularity}
s = n(\frac{1}{2} - \frac{1}{q}) -\frac{m}{r}.
\end{equation} 

This problem for $\mu =1$ has been studied by many researchers. For the Schr\"odinger operator, 
Planchon \cite{Pl} conjectured that the estimate \eqref{main_goal} is valid if and only if $r \ge 2$ and $\frac{n+1}{q} + \frac{1}{r} \le \frac{n}{2}$.
Kenig--Ponce--Vega \cite{KPV} showed the conjecture is true for $n=1$. In higher dimensions $n \ge 2$ it was proven by Vega \cite{V1} that \eqref{main_goal} holds for $q \ge \frac{2(n+2)}{n}$ and $\frac{n+1}{q} + \frac{1}{r} \le \frac{n}{2}$. When $n=2$ Rogers \cite{R} showed it for $2\le r < \infty$, $q > \frac{16}5$ and $\frac{3}{q} + \frac{1}{r} < 1$, and later the excluded endline $\frac{3}{q} + \frac{1}{r} = 1$ was obtained by Lee--Rogers--Vargas \cite{LRV}. When $n \ge 3$,  Lee--Rogers--Vargas \cite{LRV} improved the previous known result to $r \ge 2$, $q > \frac{2(n+3)}{n+1}$ and $\frac{n+1}{q} + \frac{1}{r} = \frac{n}{2}$. 
Recently it is shown by Du--Kim--Wang--Zhang \cite{DKWZ} that  the estimate \eqref{main_goal} with $r=\infty$, that is, the maximal estimate fails for $n \ge 3$. 
For a case of the wave operator it is known that  \eqref{main_goal} holds for $(r,q)$ pairs such  $2 \le r \le q$, $q \neq \infty$ and $\frac{1}{r} + \frac{n-1}{2q} \le \frac{n-1}{4}$ (see \cites{GV, KT, Pe, S}). Particularly, when $r=\infty$, Rogers--Villarroya \cite{RV} showed that \eqref{main_goal} with regularity $s> n(\frac{1}{2} - \frac{1}{q})-\frac{1}{r}$ is valid for $q \ge \frac{2(n+1)}{n-1}$.
For the fractional Schr\"odinger operator the known range of $(r,q)$ for which the estimates hold is that $2 \le r \le q$, $q \neq \infty$ and $\frac{n}{2q} + \frac{1}{r} \le \frac{n}{4}$ (see \cites{C,T2,Pa,CHKL,CO}).

The case of $\mu > 1$ has an interesting in its own right. The solution $u$ is formally written as
    \begin{equation*}
    u(t,x) = e^{it\Phi(D)}f(x) :=  \frac{1}{(2\pi)^n}\int_{\mathbb R^n} e^{i(x \cdot \xi + t \Phi(\xi))} \hat f(\xi) d\xi.
    \end{equation*}
From this form  we see that the space-time Fourier transform of $u$ is supported in the surface $S=\{(\xi, \Phi(\xi))\}$. It is known that the operator $u$ is related to the curvature of $S$ such as the sign of Gaussian curvature and the number of nonvanishing principle curvature. The Schr\"odinger operator corresponds to a paraboloid which has a positive Gaussian curvature, and the wave operator corresponds to a cone whose Gaussian curvature is zero. We are also interested in operators corresponding to a surface with negative Gaussian curvature. When $\mu >1$ there is a surface with negative Gaussian curvature. For instance, the surface $\{(\xi_1,\xi_2,\xi_1^4+2\xi_1^3\xi_2 - 2\xi_1\xi_2^3 +\xi_2^4) \}$ has negative Gaussian curvature on a neighborhood of the point $(1,0,1)$. 

In this paper we will establish a local-to-global approach as follows.

\begin{thm}\label{main_thm2}
Let $\mathbb I = (0,1)$ be a unit interval and $\mathbb B = B(0,1)$ a unit ball in $\mathbb R^n$.
Let $q_0, r_0 \in [2,\infty)$, $s(r,q)$ defined as \eqref{regularity} and $\Phi$ satisfy Condition \ref{homog_cond}.
Suppose that the local estimate
	\begin{equation}\label{thm_assum}
	\|e^{it\Phi(D)} f \|_{L_x^{q_0}(\mathbb{B};L_t^{r_0}(\mathbb{I}))}
	\leq C_{\epsilon} \|f\|_{H^{s(r_0,q_0)+\epsilon}(\mathbb{R}^n)}
	\end{equation}
holds for all $\epsilon>0$.
Then for any $q> q_0$ and $r> r_0$, the global estimate
	\begin{equation}\label{thm_goal}
	\|e^{it\Phi(D)} f \|_{L_x^q(\mathbb R^n;L_t^r(\mathbb R))}
	\le C \|f\|_{\dot H^{s(r,q)}(\mathbb R^n)}
	\end{equation}
holds, where $H^s(\mathbb R^n)$ denotes the inhomogeneous $L^2$-Sobolev space of order $s$ and $\dot H^s(\mathbb R^n)$ denotes homogeneous one.
\end{thm}

The maximal estimate, which is \eqref{thm_assum} with $r_0=\infty$, is related to pointwise convergence problems.
When $n=2$ it was proven that the maximal estimates with $m>1$ and $\mu=1$ are valid for $q_0=3$ and $s>\frac{1}{3}$ (see \cites{CK, DGL}). By interpolating with a Strichartz estimate
	\begin{equation*}
	\|e^{it\Phi(D)} f \|_{L_x^6(\mathbb{R}^2;L_t^2(\mathbb{R}))}
    \le \|e^{it\Phi(D)} f \|_{L_t^2(\mathbb{R};L_x^6(\mathbb{R}^2))}
	\leq C\|f\|_{\dot{H}^{(2-m)/4}(\mathbb{R}^2)},
	\end{equation*}
we have \eqref{thm_goal} for the line $\frac{3}{q} + \frac{1}{r} = 1$ with $r \ge 2$.
By Theorem \ref{main_thm2}, we can obtain the following global space-time estimates which is the Planchon conjecture  for $n=2$ except the endline.

\begin{cor} \label{main_thm}
Let $m>1$ and $\mu=1$. For $2 \leq r < \infty$ and $\frac{3}{q} + \frac{1}{r} < 1$, the global estimate
    \begin{equation*} 
    \|e^{it\Phi(D)} f \|_{L_x^q(\mathbb R^2;L_t^r(\mathbb R))}
    \le C \|f\|_{\dot H^{1 - \frac{2}{q} -\frac{m}{r}}(\mathbb R^2)}.
    \end{equation*}
\end{cor}

\textit{Notation.} Throughout this paper let $C>0$ denote various constants that vary from line to line, which possibly depend on $n$, $q$, $r$, $m$ and $\mu$. We use $A \lesssim B$ to denote $A \le CB$, and if $A \lesssim B$ and $B \lesssim A$ we denote by $A \sim B$.

\section{Proof of Theorem \ref{main_thm2}}
In this section we prove Theorem \ref{main_thm2} by using two propositions. In subsection 2.1 we consider a Littlewood--Paley type inequality by which the initial data $f$ can be assumed to be Fourier supported in $\{1/2 \le |\xi| \le 2\}$. In subsection 2.2 we prove a mixed norm version of Tao's $\varepsilon$-removable lemma by which the global estimates with a compact Fourier support are reduced to a local ones. In subsection 2.3 we show the two propositions imply Theorem \ref{main_thm2}.

\subsection{A Littlewood-Paley type inequality}
We discuss a Littlewood-Paley type inequality  for the operator $e^{it\Phi(D)}$ in a time variable.

Let a cut-off function $\phi \in C_0^{\infty} \big( [\frac{1}{2}, 2] \big)$ satisfy
$ \sum_{k\in\mathbb{Z}} \phi( {2^{-k}} x) =1 $.
We define Littlewood-Paley projection operators $P_{k}$ and $\widetilde{P_k}$ by
    \[
    \widehat{P_k f}(\xi) = \phi(2^{-k}|\xi|) \hat{f}(\xi)
    \quad\mbox{and}\quad
    \widehat{\widetilde{P_k} f}(\tau) = \phi (2^{-mk}|\tau| ) \hat{f}(\tau)
    \] for $\xi \in \mathbb R^n$ and $\tau \in \mathbb R$, respectively.

\begin{lem}\label{LP_lemma}
Suppose that $\Phi$ satisfies Condition \ref{homog_cond}. Then for $1<r<\infty$,
    \begin{equation*}
    \big\|e^{it\Phi(D)} f(x) \big\|_{L_t^r(\mathbb R)}
    \leq C_{m,\mu} \Big\| \Big(\sum_{\substack{j,k \in \mathbb Z: \\|k-j| \leq \frac{\log_2 \mu}{m}+2}} | \widetilde{P_j} e^{it\Phi(D)}P_kf(x) |^2\Big)^{1/2} \Big\|_{L_t^r(\mathbb R)}
    \end{equation*}
for all functions $f$ and all $x \in \mathbb R^n$.
\end{lem}

\begin{proof}
For simplicity,
    \[
    F (t) := e^{it\Phi(D)} f(x) \quad\mbox{and}\quad F_k (t) := e^{it\Phi(D)}P_kf(x).
    \]
Since the projection operators are linear, we have an identity
    \[
    F(t) = \sum_{j\in \mathbb{Z}} \widetilde{P_j}F(t) = \sum_{j\in \mathbb{Z}} \sum_{k\in \mathbb{Z}} \widetilde{P_j}F_k(t).
    \]
For any test function $\psi \in C_0^{\infty} \big( [-2, 2] \big)$ with $\psi = 1$ in $[-1,1]$,
the Fourier transform $\widehat{f}$ of $f$ is defined by
    \[
    \widehat{f}(\tau) = \lim_{R\rightarrow\infty} \frac{1}{2\pi} \int_{\mathbb{R}} e^{it\tau} \psi\Big(\frac{t}{R}\Big) f(t) dt
    \]
in the distributional sense.

We claim that $\widetilde{P_j}F_k(t) =0$
if
    \begin{equation}\label{freq_cond_null}
    |k -j| >  \frac{\log_2 \mu}{m}+2.
    \end{equation}
   Indeed, using the above definition of the Fourier transform we can write
    \[
    \begin{aligned}
    \widehat{\widetilde{P_j}F_k}(\tau)
    &= \frac{1}{(2\pi)^{n+1}} {\phi}\Big(\frac{|\tau|}{2^{mj}} \Big)  \lim_{R\rightarrow\infty} \int_{\mathbb{R}^n} e^{ix\cdot\xi} \bigg( \int_{\mathbb{R}} e^{it\tau} e^{it\Phi(\xi)} \psi\Big(\frac{t}{R}\Big) dt \bigg) \phi\Big(\frac{|\xi|}{2^k}\Big) \hat{f}(\xi) d\xi .
    \end{aligned}
    \]
In the right side of the above equation, we see that the range of $(\tau, \xi)$ is contained in 
    \[
    2^{m(j-1)}\leq |\tau|\leq 2^{m(j+1)} \quad\mbox{and}\quad 2^{(k-1)}\leq |\xi| \leq 2^{(k+1)}.
    \]
From Condition \ref{homog_cond}
we have a bound
    \[
    {\mu^{-1}} 2^{m(k-1)}\leq |\Phi(\xi)| \leq \mu 2^{m(k+1)}.
    \]
Then it follows that for  $k$ and $j$ satisfying \eqref{freq_cond_null},
\[
|\tau + \Phi(\xi)| >0.
\] 
By the integration by parts it implies that there exists a constant $C_0>0$ such that
    \[
    \Big| \int_{\mathbb{R}} e^{it\tau} e^{it\Phi(\xi)} \psi\Big(\frac{t}{R}\Big) dt \Big| \le \frac{1}{C_0 R}.
    \]
From this estimate and the Lebesgue dominated convergence theorem we obtain 
$\widehat{\widetilde{P_j}F_k}=0,$
which implies the claim.
%

By the claim, the Littlewood-Paley theory and the Cauchy-Schwarz inequality,
    \begin{equation*}
    \begin{aligned}
    \big\|e^{it\Phi(D)} f(x) \big\|_{L_t^r(\mathbb R)}
    &= \Big\| \sum_{j \in \mathbb{Z}} \widetilde{P_j} \Big( \sum_{k \in \mathbb Z} F_k(\cdot,x) \Big) \Big\|_{L_t^r(\mathbb R)} \\
    &\leq C \Big\| \Big(\sum_{j \in \mathbb{Z}} \Big| \sum_{k \in \mathbb Z: |k -j| \le  \frac{\log_2 \mu}{m}+2} \widetilde{P_j} F_k(\cdot,x) \Big|^2\Big)^{1/2} \Big\|_{L_t^r(\mathbb R)} \\
    &\leq C_{m,\mu} \Big\| \Big( \sum_{j \in \mathbb{Z}} \sum_{k \in \mathbb Z: |k -j| \le  \frac{\log_2 \mu}{m}+2} | \widetilde{P_j} F_k(\cdot,x)|^2\Big)^{1/2} \Big\|_{L_t^r(\mathbb R)}.
    \end{aligned}
    \end{equation*}
This is the desired inequality.
\end{proof}


Using the above lemma we can have the following proposition.
\begin{prop}\label{LP_prop}
Let $2 \leq q,r < \infty$. Suppose that $\Phi$ satisfies Condition \ref{homog_cond}.
If the estimate
    \begin{equation}\label{LP_prop_supp}
    \|e^{it\Phi(D)} f \|_{L_x^q(\mathbb R^n;L_t^r(\mathbb R))}
    \leq C \| f\|_{L^2(\mathbb R^n)}
    \end{equation}
holds for all $f$ with $\supp \hat f \subset \{1/2 \le |\xi| \le 2 \}$, then the estimate
    $$
    \|e^{it\Phi(D)} f \|_{L_x^q(\mathbb R^n;L_t^r(\mathbb R))}
    \leq C_{m,\mu} \| f\|_{\dot{H}^{\frac{n}{2}-\frac{n}{q}-\frac{m}{r}}(\mathbb R^n)}
    $$
holds for all $f$.
\end{prop}

\begin{proof}
The Minkowski inequality and Lemma \ref{LP_lemma} allow that
    \[
    \big\|e^{it\Phi(D)} f \big\|_{L_x^q(\mathbb R^n;L_t^r(\mathbb R))}
    \leq C_{m,\mu} \bigg\|\bigg(  \sum_{|k -j| \le  \frac{\log_2 \mu}{m}+2} \Big\| \widetilde{P_j} \big( e^{it\Phi(D)} P_kf \big) \Big\|_{L_t^r(\mathbb R)}^2  \bigg)^{1/2} \bigg\|_{L^q_x(\mathbb R^n)}.
    \]
Since $\widetilde{P_j}$ is bounded in $L^p$, it is bounded by
    \[
    C_{m,\mu} \bigg\|\bigg(  \sum_{k \in \mathbb Z} \big\| e^{it\Phi(D)} P_kf  \big\|_{L_t^r(\mathbb R)}^2  \bigg)^{1/2} \bigg\|_{L^q_x(\mathbb R^n)}.
    \]
By the Minkowski inequality we thus have
    \[
    \big\|e^{it\Phi(D)} f \big\|_{L_x^q(\mathbb R^n;L_t^r(\mathbb R))}
    \leq C_{m,\mu} \bigg( \sum_{k\in \mathbb{Z}}  \big\| e^{it\Phi(D)} P_kf \big\|_{L_x^q(\mathbb R^n;L_t^r(\mathbb R))}^2 \bigg)^{1/2}.
\]
Apply \eqref{LP_prop_supp} to the right side of the above estimate after parabolic rescaling. Then we obtain
    \begin{align*}
    \|e^{it\Phi(D)} f \|_{L_x^q(\mathbb R^n;L_t^r(\mathbb R))}
    &\leq C_{m,\mu} \bigg(\sum_{k\in\mathbb{Z}} 2^{2k(\frac{n}{2}-\frac{n}{q}-\frac{m}{r} )}\| P_kf \|_2^2 \bigg)^{1/2} \\
    &= C_{m,\mu} \| f\|_{\dot{H}^{\frac{n}{2}-\frac{n}{q}-\frac{m}{r}}(\mathbb R^n)}.
    \end{align*}
\end{proof}

\subsection{Local-to-global arguments}
We will show that the global estimate  
\eqref{LP_prop_supp} is obtained from its local estimate. Adopting the arguments in \cite{T1}, we consider the dual estimate of \eqref{LP_prop_supp}.

Let $S=\{(\xi, \Phi(\xi)) \in \mathbb R^n \times \mathbb R: 1/2 \le |x| \le 2 \}$ be a compact hypersurface with the induced (singular) Lebesgue measure $d\sigma$.
We define the Fourier restriction operator $\mathfrak R$ for a compact surface $S$ by the restriction of $\hat f$ to $S$, i.e.,
\[
\mathfrak Rf = \hat f \big|_S.
\]
Its adjoint operator $\mathfrak R^*f = \widehat{fd\sigma}$ can be viewed as $e^{it\Phi(D)}\hat{g}$, where the Fourier transform $\hat g(\xi)$ of $g$ corresponds to $f(\xi, \Phi(\xi))$.

%
Let $\rho>0$ be the decay of $\widehat{d\sigma}$, i.e.,
\begin{equation} \label{surdecay}
|\widehat{d\sigma}(x)| \lesssim (1+|x|)^{-\rho}, \qquad x \in \mathbb R^{n+1}.
\end{equation}
It is known that $\rho$ is determined by the number of nonzero principal curvatures of the surface $S$, which is equal to the rank of the Hessian $H_{\Phi}$. Specifically, if $H_{\Phi}$ has rank at least $k$ then
\[
\rho = k/2,
\]
see \cite{St}*{subsection 5.8, VIII}. From Condition \ref{homog_cond} we have $k \ge 1$.

When a function $f$ has a compact Fourier support, the $\widehat{fd\sigma}$ decays away from the support of $\hat f$ because of the decay of $\widehat{d\sigma}$. Thus if $f$ and $g$ are compactly Fourier supported and their supports are far away from each other then the interaction between $\widehat{fd\sigma}$ and $\widehat{gd\sigma}$ is negligible. 
\begin{defn}
A finite collection $\{Q(z_i,R)\}_{i=1}^{N}$ of balls in $\mathbb{R}^{n+1}$ with radius $R>0$ is called  $(N,R)$-\textit{sparse} if the centers $\{z_i\}$ are $(NR)^\gamma$-separated where $\gamma := n/\rho~ (\ge 2)$. 
\end{defn}

From the definition of $(N,R)$-sparse we have a kind of orthogonality as follows. Let $\phi$ be a radial Schwartz function which is positive on the ball $B(0,3/2)$ and $\phi = 1$ on the unit ball $B(0,1)$ and whose Fourier transform is supported on the ball $B(0,2/3)$. 

\begin{lem}[\cite{T1}*{in the proof of Lemma 3.2}]  \label{lem:spase_decp}
Let $\{Q(z_i,R)\}_{i=1}^{N}$ be a $(N,R)$-sparse collection 
and $\phi_i(z)=\phi(R^{-1}(z-z_i))$ for $i=1,\cdots, N$. Then there is a constnat $C$ independent of $N$ such that
\begin{equation} \label{eqn:dep}
\Big\| \sum_{i=1}^{N} f_i \ast \hat \phi_i \big|_S  \Big\|_2 \le CR^{1/2} \Big( \sum_{i=1}^{N} \|f_i\|_2^2 \Big)^{1/2}
\end{equation}
for all $f_i \in L^2(\mathbb R^{n+1})$.
\end{lem}
A proof of the above lemma is given in Appendix. 

Let $\mathbb I_R=(0,R)$ denote an $R$-interval and $\mathbb B_R$ the ball of radius $R$ centered at the origin in $\mathbb R^n$. Using Lemma \ref{lem:spase_decp} we have an intermediate result.
\begin{prop} \label{lem:sparseEst}
Let $R>0$ and $1 < q,r \leq 2$.
Suppose that there is a constant $A(R)$ such that
    \begin{equation}\label{eqn:loc_rest}
    \|\mathfrak R(\chi_{\mathbb{I}_R \times \mathbb B_R} f) \|_{L^2(d\sigma)}
    \leq A(R) \|f\|_{L_{x}^{q}(\mathbb R^n;L_{t}^{ r}(\mathbb R))}
    \end{equation}
for all $f \in L_{x}^{q}(\mathbb R^n;L_{t}^{r}(\mathbb R))$.
Then for any $(N,R)$-sparse collection $\{Q(z_i,R)\}_{i=1}^{N}$ there is a constant $C$ independent of $N$ such that
    \begin{equation}\label{eqn:loc_rest_result}
    \|\mathfrak Rf\|_{L^2(d\sigma)}
    \leq C A(R) \|f\|_{L_{x}^{q}(\mathbb R^n;L_{t}^{ r}(\mathbb R))}
    \end{equation}
for all $f$ supported in $\cup_{i=1}^{N} Q(z_i,R)$.
\end{prop}

\begin{proof}
Let $f_i = f \chi_{Q(z_i,R)}$. Then,
\[\mathfrak R f_i = \hat f_i \big|_S =\widehat{ f_i \phi_i}\big|_S =  (\hat f_i \ast \hat \phi_i) \big|_S,\] where $\phi_i(z)$ is defined as in Lemma \ref{lem:spase_decp}. Since $\hat \phi_i$ is supported on the ball $B(0,\frac{2}{3R})$, we may restrict the support of $\hat f_i$ to a $O(1/R)$-neighborhood of the surface $S$ and write
\[
\mathfrak R f_i  =  (\hat f_i \big|_{\mathcal N_{1/R} (S)} \ast \hat \phi_i) \big|_S
\]
where $\mathcal N_{1/R} (S)$ is a $O(1/R)$-neighborhood of the surface $S$. Let $\tilde{\mathfrak R}$ be another restriction operator defined by $\tilde{\mathfrak R} f = \hat f \big|_{\mathcal N_{1/R}(S)}$. If $f$ is suppported in $\cup_{i=1}^{N} Q(z_i,R)$, we write
\[
\mathfrak R f = \sum_{i=1}^{N} (\tilde{\mathfrak R}f_i \ast \hat \phi_i) \big|_{S}.
\]
By Lemma \ref{lem:spase_decp},
\[
\|\mathfrak R f \|_{L^2(d\sigma)} \le C R^{1/2} \Big( \sum_{i=1}^{N} \| \tilde{\mathfrak R} f_i \|_{L^2(\mathcal N_{1/R} (S))}^2 \Big)^{1/2}.
\]
Since the estimate \eqref{eqn:loc_rest} is translation invariant, by a slice argument we have
\[
\| \tilde{\mathfrak R} f_i \|_{L^2(\mathcal N_{1/R} (S))}
\le CR^{-1/2} A(R)  \|f_i\|_{L_x^{q}(\mathbb R^n;L_t^r(\mathbb R))}.
\]
By combining the previous two estimates,
\[
\|\mathfrak R f \|_{L^2(d\sigma)} \le C A(R)  \Big( \sum_{i=1}^{N} \|f_i \|^2_{L_{x}^{q}(\mathbb R^n;L_{t}^{r}(\mathbb R))} \Big)^{1/2}.
\]
If $1 \le r \le q \le 2$ then by $\ell^r \subset \ell^{q} \subset \ell^{2}$, 
\begin{align*}
\Big( \sum_{i=1}^{N} \|f_i \|^2_{L_{x}^{q}(\mathbb R^n;L_{t}^{r}(\mathbb R))} \Big)^{1/2}
&\le \Big( \sum_{i=1}^{N} \|f_i \|^q_{L_{x}^{q}(\mathbb R^n;L_{t}^{r}(\mathbb R))} \Big)^{1/q} \\
&=\Big(\int_{\mathbb R^n}\sum_{i=1}^{N}\|f_i\|_{L_t^{r}(\mathbb R)}^{q} dx\Big)^{1/q} \\
&\le \Big(\int_{\mathbb R^n} \Big(\sum_{i=1}^{N}\|f_i\|_{L_t^{r}(\mathbb R)}^{r} \Big)^{q/r} dx \Big)^{1/q} \\
&=\|f\|_{L_{x}^{q}(\mathbb R^n;L_{t}^{ r}(\mathbb R))}.
\end{align*}
If $1 \le q \le r \le 2$ one can use the embedding $\ell^r \subset \ell^2$ and the Minkowski inequality to get
\begin{align*}
\Big( \sum_{i=1}^{N} \|f_i \|^2_{L_{x}^{q}(\mathbb R^n;L_{t}^{r}(\mathbb R))} \Big)^{1/2}
&\le \Big( \sum_{i=1}^{N} \|f_i \|^r_{L_{x}^{q}(\mathbb R^n;L_{t}^{r}(\mathbb R))} \Big)^{1/r} \\
&\le \Big(\int_{\mathbb R^n} \Big(\sum_{i=1}^{N}\|f_i\|_{L_t^{r}(\mathbb R)}^{r} \Big)^{q/r} dx\Big)^{1/q} \\
&=\|f\|_{L_{x}^{q}(\mathbb R^n;L_{t}^{ r}(\mathbb R))}.
\end{align*}
Therefore we have \eqref{eqn:loc_rest_result}.
\end{proof}

We now extend the $(N,R)$-sparse sets to the whole space. For this we need the following decomposition lemma.
\begin{lem}[\cite{T1}]\label{lem:DecE}
Let $E$ be a subset in $\mathbb R^n$ with $|E| >1$. Suppose that $E$ is a finite union of finitely overlapping cubes of side-length $c \sim 1$.
Then for each $K \in \mathbb{N}$, there are subsets $E_1, E_2, \cdots, E_K$ of $E$ with
\[
    E = \bigcup_{k=1}^{K} E_k
\]
such that each $E_k$ has $O(|E|^{1/K})$ number of  $(O(|E|), |E|^{O(\gamma^{k-1})})$-sparse collections 
\[
\mathbf S_1, \mathbf S_2, \cdots, \mathbf S_{O(|E|^{1/K})}
\] of which the union $\mathbf S_1 \cup \mathbf S_2 \cup  \cdots \cup \mathbf S_{O(|E|^{1/K})}$ is a covering of $E_k$.
\end{lem}
This lemma is a precise version of Lemma 3.3 in \cite{T1}. A detailed proof can be found in Appendix.

Using the above lemma we have the following proposition.
\begin{prop}\label{prop:glob_d}
Let $1 < q_0,r_0 < \infty$.
Suppose that for any $\epsilon>0$ and any $(N,R)$-sparse collection $\{Q(z_i,R)\}_{i=1}^{N}$ in $\mathbb R^{n+1}$,
the estimate
    \begin{equation} \label{eqn:sparse_rest}
    \| \mathfrak R f \|_ {L^{2}(d\sigma)}\leq C_{\epsilon} R^{\epsilon} \|f\|_{L_x^{q_0}(\mathbb R^n; L_t^{r_0}(\mathbb R) )}
    \end{equation}
holds for all $f$ supported in $\cup_{i=1}^{N} Q(z_i,R)$.
Then for any $1\leq q < q_0$ and $1\leq r < r_0$, the estimate
    \begin{equation*} 
    \| \mathfrak R f \|_ {L^{2}(d\sigma)}\le C\|f\|_{L_x^{q}(\mathbb R^n; L_t^{r}(\mathbb R) )}
    \end{equation*}
holds for all $f \in L_x^{q}(\mathbb R^n; L_t^{r}(\mathbb R) )$.
\end{prop}

\begin{proof}
By interpolation (see \cite{F}),
it suffices to show that for $1 \leq q < q_0$ and $1 \leq r < r_0$,
the restricted type estimate
    \begin{equation*} 
    \|\mathfrak R \chi_E \|_{L^2(d\sigma)}
    \leq C   \|\chi_E\|_{L^{q}(\mathbb R^n;L^{r}(\mathbb R))}
    \end{equation*}
for all subset $E$ in $\mathbb R^{n+1}$. We may assume $|E| > 1$, otherwise the estimate is trivial.  
Since the set $S$ is compact,
$\chi_E$ can be replaced with $\chi_E \ast \varphi$, where $\varphi$ is a bump function supported on a cube of sidelength $c \sim 1$ such that $\hat \varphi$ is positive on $S$. Thus, we may further assume that $E$ is the union of $c$-cubes. 

We denote by $\mathrm{proj} (E)$ the projection of $E$ onto the $x$-plane. For each grid point $x \in c\, \mathbb Z^n \cap \mathrm{proj} (E)$, we define $E_x$ to be the union of $c$-cubes in $E$ that intersect $\mathbb R \times \{x\}$.  Let $E^j$ be the union of $E_x$ which satisfies 
           \[2^{j-1} <  \text{ the number of}\,\, c\,\text{- cubes contained in}\,\, E_x \le 2^{j+1} \]
for $j \in \mathbb N$, (see Figure \ref{fig1}).
Then,
    $$
    E= \bigcup_{j \ge 1} E^j.
    $$

\begin{figure} 
 \begin{center} 
  {\includegraphics[width=0.9\textwidth]{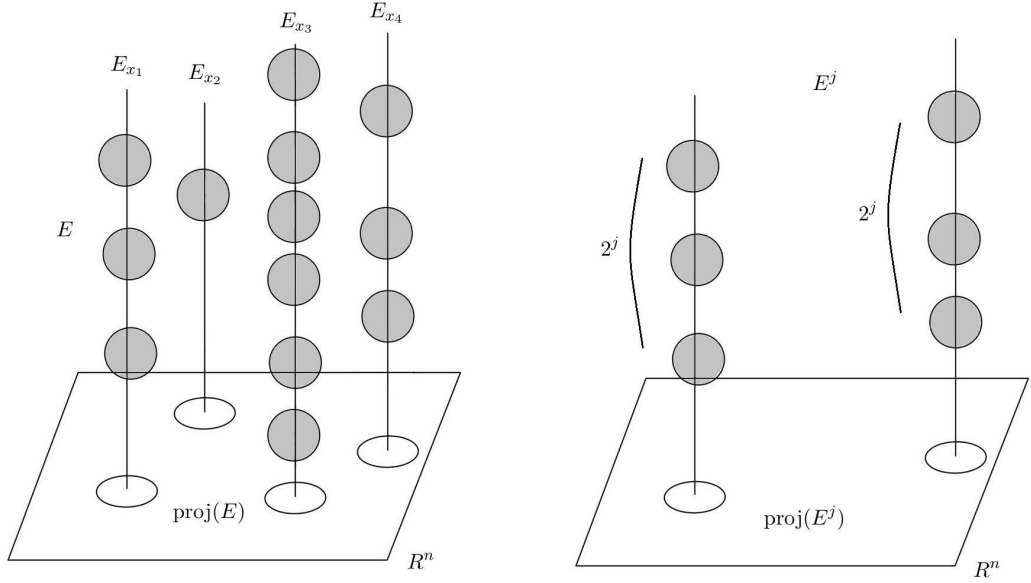}}
 \end{center}
\caption{The sets $E$, proj$E$, $E_x$ and $E^j$ in the proof of Proposition \ref{prop:glob_d}.}
\label{fig1}
\end{figure}

By using Lemma \ref{lem:DecE} with 
\[
K := \frac{\log(1/\epsilon)}{2\log \gamma}  + 1,
\] 
the $E^j$ is decomposed into $E^j_k$'s which are covered by $O(|E^j|^{1/K})$ number of $(O(|E^j|),|E^j|^{C\gamma^{k-1}}))$-sparse collections. We apply \eqref{eqn:sparse_rest} to these sparse collections and obtain
    $$
    \|\mathfrak R \chi_{E^j_k} \|_{L^2(d\sigma)}
    \leq C_{\epsilon} |E^j|^{1/K} (|E^j|^{C\gamma^{k-1}})^{\epsilon}  \|\chi_{E^j_k} \|_{L_x^{q_0}(\mathbb R^n ; L_t^{r_0}(\mathbb R))}.
    $$
Summing over k, we have
    \begin{align*}
    \|\mathfrak R \chi_{E^j} \|_{L^2(d\sigma)}
    &\leq \sum_{k=1}^K \|\mathfrak R \chi_{E^j_k} \|_{L^2(d\sigma)} \\
    &\leq C_{\epsilon}  |E^j|^{1/K} (|E^j|^{C\gamma^{K-1}})^{\epsilon} \|\chi_{E^j} \|_{L_x^{q_0}(\mathbb R^n ; L_t^{r_0}(\mathbb R))}
    \end{align*}
where $K$ is absorbed into $C_\epsilon$.
Since $|E^j| \leq 2^{j+1} |\mathrm{proj\,} (E^j)|$, we have
    \[
    \|\mathfrak R \chi_{E^j} \|_{L^2(d\sigma)}
    \le C_{\epsilon} 2^{j(\frac{1}{r_0} + \delta(\epsilon) )} |\mathrm{proj} (E^j)|^{\frac{1}{q_0} + \delta(\epsilon)},
    \]
where  
\[
\delta(\epsilon) :=\frac{1}{K} +C\gamma^{K-1}\epsilon.
\]
Since $\lim_{\epsilon \to 0}\delta(\epsilon) = 0$,
we can take $\epsilon >0$ such that
\[
0< \delta(\epsilon)+\epsilon \le \min \bigg(\frac{1}{q}-\frac{1}{q_0}, \frac{1}{r} - \frac{1}{r_0} \bigg).
\]
Thus,
    \begin{align*}
    \|\mathfrak R \chi_E \|_{L^2(d\sigma)}
    &\leq \sum_{j\geq1} \|\mathfrak R \chi_{E^j} \|_{L^2(d\sigma)} \\
    &\leq C_{\epsilon} \sum_{j\geq1} 2^{j(\frac{1}{r_0} + \delta(\epsilon))} |\mathrm{proj} (E^j)|^{\frac{1}{q_0} + \delta(\epsilon)}\\
&\leq C  \sum_{j\geq1} 2^{-\epsilon j} 2^{\frac{1}{r} j} |\mathrm{proj} (E^j)|^{\frac{1}{q}} \\
&\leq C  \sum_{j\geq1} 2^{-\epsilon j} \|\chi_{E}\|_{L_x^{q}(\mathbb R^n;L_t^{r}(\mathbb R))}  \\
 &\leq C \|\chi_{E}\|_{L_x^{q}(\mathbb R^n;L_t^{r}(\mathbb R))}.
\end{align*}
\end{proof}

Combining Proposition \ref{lem:sparseEst} and Proposition \ref{prop:glob_d} we obtain an extension of Tao's epsilon removal lemma as follows.
\begin{prop}\label{eps_remv}
Let $1 < q_0,r_0 \leq 2$.
Suppose that
    $$
    \|\mathfrak R(\chi_{\mathbb{I}_R \times \mathbb B_R} f) \|_{L^2(d\sigma)}
    \leq C_{\epsilon} R^{\epsilon} \|f\|_{L_{x}^{q_0}(\mathbb R^n;L_{t}^{r_0}(\mathbb R))}
    $$
for all $\epsilon>0$, $R>1$ and all $f \in L_{x}^{q}(\mathbb R^n;L_{t}^{ r}(\mathbb R))$.
Then for any $1\leq q < q_0$ and $1\leq r < r_0$, 
    $$
    \| \mathfrak R f \|_ {L^{2}(d\sigma)}\le C\|f\|_{L_x^{q}(\mathbb R^n; L_t^{r}(\mathbb R) )}
    $$
for all $f \in L_x^{q}(\mathbb R^n; L_t^{r}(\mathbb R) )$.
\end{prop}

%

Now we are ready to prove Theorem \ref{main_thm2}. The theorem follows from Proposition \ref{LP_prop} and Proposition \ref{eps_remv} as follows.

\subsection{Proof of Theorem \ref{main_thm2}}

Let $P_0$ be the Littlewood-Paley projection operator as in subsection 2.1.
By rescaling $x \mapsto 2^{-k}x$ and $t \mapsto 2^{-mk}t$,  the estimate \eqref{thm_assum} implies
    $$
    \| e^{it\Phi(D)} P_0 f\|_{L_x^{q_0}(\mathbb{B}_{R}; L_t^{r_0}(\mathbb{I}_{R^{m}}) )}
    \leq C_{\epsilon} 2^{k\epsilon} \|P_0 f\|_{L^2(\mathbb R^n)}
    $$
for all $k\geq1$ and $\epsilon>0$. 
Since $m\geq 1$, we have
    $$
    \| e^{it\Phi(D)} P_0 f\|_{L_x^{q_0}(\mathbb{B}_{R}; L_t^{r_0}(\mathbb{I}_{R}) )}
    \leq C_{\epsilon} 2^{k\epsilon} \|P_0 f\|_{L^2(\mathbb R^n)} .
    $$
By Proposition \ref{eps_remv} and duality,
    $$
    \| e^{it\Phi(D)} P_0 f\|_{L_x^{q}(\mathbb{R}^n; L_t^{r}(\mathbb{R}) )}
    \leq C\|P_0 f\|_{L^2(\mathbb R^n)}.
    $$
By Proposition \ref{LP_prop}, we obtain the desired estimate. \qed

\section{Appendix} \label{sec:append}
\subsection{Proof of Lemma \ref{lem:spase_decp}}
We divide the left side of \eqref{eqn:dep} into two parts
\[
\| \sum_{i=1}^{N} f_i \ast \hat \varphi_i |_{S} \|_2^2 
= \sum_{i} \|f_i \ast \hat \varphi_i |_{S} \|_2^2
+ \sum_{i \neq j} \int f_i \ast \hat \varphi_i  \overline{f_j \ast \hat \varphi_j} d\sigma.
\]
By a basic restriction estimate we have $\| f_i \ast \hat \varphi_i |_{S}\|_2  \lesssim R^{1/2} \|f_i\|_2$. Thus,
\begin{equation} \label{Bres}
\sum_{i=1}^{N} \| f_i \ast \hat \varphi_i |_{S}\|_2^2  \lesssim R \sum_{i=1}^{N}  \|f_i\|_2^2.
\end{equation}
By Parseval's identity,
\[
 \int f_i \ast \hat \varphi_i  \overline{f_j \ast \hat \varphi_j} d\sigma 
= \int \overline{\check f_j \varphi_j} ((\check f_i \varphi_i)  \ast \widehat{d\sigma}),
\] where the $\check{ }$ denotes the inverse Fourier transform.
It is bounded by
\[
\big(\sup_{z,w} | \varphi_j^{1/2}(z)  \varphi^{1/2}_i(w)   \widehat{d\sigma}(z-w) | \big)  \|\check f_i \varphi^{1/2}_i\|_1 \|\check f_j \varphi^{1/2}_j\|_1.
\]
By the Cauchy-Schwarz inequality and Plancherel's theorem,
\[
\|\check f_i \varphi^{1/2}_i \|_{1} \lesssim R^{(n+1)/2} \| f_i\|_{2}.
\]
By \eqref{surdecay},
\[
\sup_{z,w} | \varphi^{1/2}_j(z)  \varphi^{1/2}_i(w)   \widehat{d\sigma}(z-w) | \lesssim |z_i-z_j-2R|^{-\rho}.
\]
Since $|z_i-z_j| \ge (NR)^{\gamma}$ and $\gamma \ge 2$, we have that $|z_i-z_j-2R|$ is comparable to $|z_i-z_j|$. Thus,
\[
\sup_{z,w} | \varphi^{1/2}_j(z)  \varphi^{1/2}_i(w)   \widehat{d\sigma}(z-w) | \lesssim |z_i-z_j|^{-\rho}.
\]
Combining these estimates we have
\begin{align*}
\sum_{i \neq j} \int f_i \ast \hat \varphi_i  \overline{f_j \ast \hat \varphi_j} d\sigma &\lesssim R^{n+1} \sum_{i=1}^{N}\sum_{j=1}^{N} |z_i-z_j|^{-\rho} \| f_i\|_{2} \| f_j\|_{2} \\
&\lesssim R^{n+1} N \max_{i,j} |z_i-z_j|^{-\rho} \sum_{i=1}^{N} \| f_i\|_{2}^2.
\end{align*}
Since $|z_i-z_j| \ge (NR)^{\gamma} \ge N^{\frac{1}{\rho}}R^{\frac{n}{\rho}}$, it follows that
\[
\sum_{i \neq j} \int f_i \ast \hat \varphi_i  \overline{f_j \ast \hat \varphi_j} d\sigma \lesssim R \sum_{i=1}^{N} \| f_i\|_{2}^2.
\]
From the above estimate and \eqref{Bres} we obtain \eqref{eqn:dep}.
\qed

\subsection{Proof of Lemma \ref{lem:DecE}}
Fix $K \in \mathbb N$. We define $R_0=1$ and $R_k$ for $k=1,2,\cdots, K$ recursively by 
\begin{equation} \label{sepB}
R_{k} = |E|^{\gamma} R_{k-1}^{\gamma}.
\end{equation}
From this definition we have $R_k = |E|^{\frac{\gamma^{k+1}-\gamma}{\gamma-1}}$. Let  $E_0 = \emptyset$. We define $E_k$ for $k=1,2,\cdots, K$ to be the set of all $x \in E \setminus \cup_{j=0,1,2,\cdots,k-1} E_j$ such that
\begin{equation} \label{Econ1}
|E \cap B(x, R_k)| \le |E|^{k/K}.
\end{equation}
Then, $E = \bigcup_{k=1}^{K} E_k$. From this construction it follows that that for $x \in E_k$, $k=2,3, \cdots, K$,
\begin{equation} \label{Econ2}
|E \cap B(x, R_{k-1})| > |E|^{(k-1)/K}.
\end{equation}

We cover $E_k$ with finitely overlapping $R_k$-balls $\mathbf C_{E_k} := \{B_i=B(x_i,R_k): x_i \in E_k\}$.
Since $E$ is a finite union of cubes of side-length $c \sim 1$, it is obvious that $\#\mathbf C_{E_k} \lesssim |E|$. For each $B_i \in \mathbf C_{E_k}$ we cover $E_k \cap B_i$ with finitely overlapping $R_{k-1}$-balls $\mathbf C_{E_k \cap B_i} :=\{B'_{ij}=B'(y_j, R_{k-1}):y_j \in E_k \cap B_i \}$, that is,
\[
E_k \cap B_i = \bigcup_{B'_{ij} \in \mathbf C_{E_k \cap B_i}} E_k \cap B_{ij}'.
\]
Since $((E  \setminus E_k) \cap B_{ij}') \subset ((E  \setminus E_k) \cap B_{i})$ for all $j$, we have
\[
(E_k \cap B_i)  \cup ((E  \setminus E_k) \cap B_{i}) \supset \bigcup_{B'_{ij} \in \mathbf C_{E_k \cap B_i}} (E_k \cap B_{ij}') \cup ((E  \setminus E_k) \cap B_{ij}'),
\]
thus
\[
E \cap B_i \supset \bigcup_{B'_{ij} \in \mathbf C_{E_k \cap B_i}} E \cap B_{ij}'.
\]
By finitely overlapping,
\[
\# \mathbf C_{E_k \cap B_i} \lesssim \max_{B'_{ij} \in \mathbf C_{E_k \cap B_i}} \frac{|E \cap B_i|}{|E \cap B_{ij}'|}.
\]
By \eqref{Econ1} and \eqref{Econ2} the above is bounded by $C|E|^{1/K}$, and we have $\#\mathbf C_{E_k \cap B_i} \le C|E|^{1/K}$ for all $i$.
Thus,
\[
E_k \subset \bigcup_{i=1}^{O(|E|)} \bigcup_{j=1}^{O(|E|^{1/K})} B_{ij}'.
\] 

We choose $O(R_k)$-separated balls $\{B'_{ij(i)}\}_{i=1}^{O(|E|)}$. Then it becomes a $(O(|E|), R_{k-1})$-sparse collection because of \eqref{sepB}. Since $R_{k-1} = |E|^{O(\gamma^{k-1})}$ and every $B_i \in \mathbf C_{E_k}$ has the covering $\mathbf C_{E_k \cap B_i}$ of cardinality $O(|E|^{1/K})$, there are $O(|E|^{1/K})$ number of $(O(|E|), |E|^{O(\gamma^{k-1})})$-sparse collections $\mathbf S_1, \mathbf S_2, \cdots, \mathbf S_{O(|E|^{1/K})}$ such that 
\[
E_k \subset \bigcup_{j=1}^{O(|E|^{1/K})} \bigcup_{B' \in \mathbf S_j} B'.
\]
\qed

\end{document}